\documentclass[12pt,reqno]{amsart}
\usepackage{amssymb,amscd,amsbsy}
\newcommand{\lb}{\linebreak}

\renewcommand{\a}{\alpha}
\renewcommand{\b}{\beta}

\newcommand{\vk}{\varkappa}

\newcommand{\s}{\sigma}

\newcommand{\f}{\varphi}

\newcommand{\G}{\Gamma}
\newcommand{\D}{\Delta}

\newcommand{\B}{{\mathcal B}}
\newcommand{\cc}{{\mathcal C}}

\newcommand{\F}{{\mathcal F}}
\newcommand{\Q}{{\mathcal Q}}
\newcommand{\h}{{\mathcal H}}

\newcommand{\p}{{\mathcal P}}
\newcommand{\cp}{{\mathcal P}}
\newcommand{\cR}{{\mathcal R}}
\newcommand{\X}{{\mathcal X}}

\newcommand{\R}{{\Bbb R}}

\newcommand{\0}{{\boldsymbol{0}}}

\newcommand{\bs}{\boldsymbol}

\newcommand{\GG}{{\boldsymbol{\varGamma}}}

\newcommand{\bS}{{\boldsymbol S}}

\newcommand{\rf}[1]{(\ref{#1})}

\newcommand{\df}{\stackrel{\mathrm{def}}{=}}

\newcommand{\supp}{\operatorname{supp}}

\newcommand{\trace}{\operatorname{trace}}
\newcommand{\rank}{\operatorname{rank}}
\newcommand{\const}{\operatorname{const}}

\newcommand{\eeq}{\end{equation}}
\newcommand{\beq}{\begin{equation}}
\newcommand{\bay}{\begin{eqnarray}}
\newcommand{\ey}{\end{eqnarray}}
\newcommand{\bey}{\begin{eqnarray*}}
\newcommand{\eey}{\end{eqnarray*}}

\newcommand{\be}{\infty}

\newcommand{\bl}{\blacksquare}

\newcommand{\Pf}{{\bf Proof. }}

\newcommand{\ov}{\overline}

\newtheorem{thm}{\hspace{\parindent}Theorem}[section]

\newtheorem{lem}[thm]{\hspace{\parindent}Lemma}

\begin{document}

\newcommand{\vse}{\vspace{.2in}}
\numberwithin{equation}{section}

\title{\bf Distorted Hankel integral operators}
\author{A.B. Aleksandrov and V.V. Peller}
\thanks{The first author is partially supported by Grant 02-01-00267
of Russian Foundation of Fundamental Studies and by Grant 326.53 of Integration.
The second author is partially supported by NSF grant DMS 0200712.}
\maketitle

\begin{abstract}
For $\a,\b>0$ and for a locally integrable function (or, more generally, a distribution)
$\f$ on $(0,\be)$, we study integral ooperators ${\frak G}^{\a,\b}_\f$ on $L^2(\R_+)$
defined by $\big({\frak G}^{\a,\b}_\f f\big)(x)=\int_{\R_+}\f\big(x^\a+y^\b\big)f(y)dy$.
We describe the bounded and compact operators ${\frak G}^{\a,\b}_\f$ and operators ${\frak G}^{\a,\b}_\f$
of Schatten--von Neumann class $\bS_p$. We also study continuity properties of the averaging projection
$\Q_{\a,\b}$  onto the operators of the form ${\frak G}^{\a,\b}_\f$. In particular, we show that if $\a\le\b$ and $\b>1$, then 
${\frak G}^{\a,\b}_\f$ is bounded on $\bS_p$ if and only if $2\b(\b+1)^{-1}<p<2\b(\b-1)^{-1}$.
\end{abstract}

\section{\bf Introduction}
\setcounter{equation}{0}

\

We are going to study a class of integral operators on $L^2(\R_+)$ that generalizes Hankel integral operators.

For a function $\f$ in $L^1(\R_+)$, the integral Hankel operator $\GG_\f$ is defined on $L^2(\R_+)$ by
$$
(\GG_\f f)(x)=\int_{\R_+}\f(x+y)f(y)dy.
$$
It is easy to see that such operators are bounded on $L^2(\R_+)$. Operators $\GG_\f$ can be bounded
in a much more general case. For $\GG_\f$ to be bounded,
$\f$ does not have to be a function, it can be a distribution. Bounded Hankel operators 
$\GG_\f$ are unitarily equivalent to Hankel operators on $\ell^2$, i.e., operators  with 
Hankel matrices of the form $\{\a_{j+k}\}_{j,k\ge0}$. These facts can be found in \cite{P}, Ch. 1, \S 8.

In this paper for $\a,\,\b>0$, we study the integral operators ${\frak G}^{\a,\b}_\f$ on $L^2(\R_+)$ defined by 
\bay
\label{dis}
\big({\frak G}^{\a,\b}_\f f\big)(x)=\int_{\R_+}\f\big(x^\a+y^\b\big)f(y)dy.
\ey
Clearly, if $\f$ is a locally integrable function on $\R_+=(0,\be)$, the right-hand side of \rf{dis} is well defined for 
smooth functions $f$ with compact support in $(0,\be)$. 
The integral on the right-hand side of \rf{dis} makes sense for distributions $\f$ on $(0,\be)$
and infinitely differentiable functions $\f$ with compact support in $\R_+$.
We say that for a distribution $\f$, the operator ${\frak G}^{\a,\b}_\f$ is bounded on $L^2(\R_+)$ 
if it extends by continuity to $L^2(\R_+)$.
Integral operators ${\frak G}^{\a,\b}_\f$ are called {\it distorted Hankel operators}.
We are going to study boundedness, compactness and
Schatten-von Neumann properties of distorted Hankel operators.

Obviously, for $\a=\b=1$, the operator ${\frak G}^{\a,\b}_\f$ coincides with the integral Hankel operator $\GG_\f$.
On the other hand, the limit case of the operators ${\frak G}^{\a,\a}_\f$ as $\a\to\be$ are the integral operators
$Q_\f$ on $L^2(\R_+)$ defined by
$$
(Q_\f f)(x)=\int_{\R_+}\f(\max\{x,y\})f(y)dy.
$$
We refer the reader to \cite{AJPR} where the operators $Q_\f$ are studied in detail. Note, however, that
properties of the operators $Q_\f$ are quite different from properties of the operators ${\frak G}^{\a,\b}_\f$.

In \S 2 we collect necessary information on Schatten--von Neumann classes, \lb weighted integral
Hankel operators, and discuss properties of the averaging projection and weighted projections onto the
integral Hankel operators. We state in \S 2 known results (Theorems A--D) that will be used in Sections 4 and 5. 

In \S 3 we show that Theorem B, that describes the weighted integral Hankel operators $\GG_\f^{\a,\b}$
of class $\bS_p$, does not extend to $\a$ and $\b$ not satisfying the hypotheses of Theorem B.

In \S 4 we reduce the study of the operators ${\frak G}^{\a,\b}_\f$ to the study of weighted integral Hankel operators
(see the definition in \S 2). We describe the operators ${\frak G}^{\a,\b}_\f$ that belong to the
Schatten--von Neumann class $\bS_p$ under a certain condition on $\a,\,\b$, and $\p$.

In \S 5 we introduce the averaging projection $\Q_{\a,\b}$
onto the subspace of operators of the form ${\frak G}^{\a,\b}_\f$
and we study their metric properties. In particular, we show that if $\a\le\b$, $\b>1$, and 
$$
\frac{2\b}{\b+1}<p<\frac{2\b}{\b-1},
$$
then $\Q_{\a,\b}$ is bounded on $\bS_p$. Moreover, this result is sharp. It is interesting to observe that
both $\frac{2\b}{\b+1}$ and $\frac{2\b}{\b-1}$ go to 2 as $\b$ tends to $\be$. However, in the limiting case
the averaging projection on the subspace of the operators of the from $Q_\f$ is bounded for all $p\in(1,\be)$,
see \cite{AJPR}.

It is also interesting that, unlike in the case of Hankel integral operators, 
for $\a,\,\b\in(0,1)$, the averaging projection $\Q_{\a,\b}$ is bounded on $\bS_1$, and on the space of
bounded and compact operators.

\pagebreak

\section{\bf Preliminaries}
\setcounter{equation}{0}

\

{\bf Schatten--von Neumann classes.} Recall that for a bounded operator $T$ on a Hilbert space $\h$ the singular values
$s_j(T)$, $j\ge0$, are defined by
$$
s_j(T)=\inf\{\|T-K\|:~K\in\B(\h),~\rank K\le j\}.
$$
Here $\B(\h)$ denotes the space of bounded linear operators on $\h$.

The Schatten--von Neumann class $\bS_p=\bS_p(\h)$, $0<p<\be$, consists of the operators $T$ on $\h$ such that 
$$
\|T\|_{\bS_p}=\left(\sum_{j\ge0}\big(s_j(T)\big)^p\right)^{1/p}<\be.
$$
We denote by $\bs{\cc}(\h)$ the space of compact operators on $\h$.

If $1\le p<\be$, then $\|\cdot\|_{\bS_p}$ is a norm, which makes $\bS_p$ a Banach space. For $p<1$,
$\|\cdot\|_{\bS_p}$ does not satisfy the triangle inequality, it is a quasinorm (i.e., 
$\|T_1+T_2\|_{\bS_p}\le\const(\|T_1\|_{\bS_p}+\|T_2\|_{\bS_p})$  for $T_1,\,T_2\in\bS_p$), which makes $\bS_p$
a quasi-Banach space.

The linear functional $\trace$ is defined on $\bS_1$ by
$$
\trace T=\sum_{j\ge0}(Te_j,e_j),\quad T\in\bS_1,
$$
where $\{e_j\}_{j\ge0}$ is an orthonormal basis in $\h$. Moreover, the right-hand side does not depend on the choice of the
basis.

If $1<p<\be$, the dual space $\bS_p^*$ can be identified with $\bS_{p'}$ with respect to the pairing
\bay
\label{sled}
\langle T,R\rangle=\trace TR^*, \quad T\in\bS_p,\quad R\in\bS_{p'}.
\ey
Here $p'=p/(p-1)$ is the dual exponent. With respect to the same pairing \rf{sled} one can identify $\bS_1^*$ 
with $\B(\h)$ and $\big(\bs{\cc}(\h)\big)^*$ with $\bS_1$. 

In the case $\h=L^2(\X,\mu)$ and $\mu$ is a $\s$-finite measure, 
the space $\bS_2$ coincides with the set of integral operators ${\frak I}_k$,
$$
({\frak I}_kf)(x)=\int_{\X}k(x,y)f(y)\,d\mu(y),\quad f\in L^2(\X,\mu),\quad x,\,y\in\X,
$$
with $k\in L^2(\X\times\X,\mu\otimes\mu)$ and
$$
\|{\frak I}_k\|_{\bS_2}=\left(\,\,\iint\limits_{\X\times\X}|k(x,y)|^2d\mu(x)\,d\mu(y)\right)^{1/2}.
$$
Moreover, if $T$ and $R$ are integral operators in $\bS_2$ with kernel functions $k$ and $\xi$, then
$$
\trace TR^*=\iint\limits_{\X\times\X}k(x,y)\ov{\xi(x,y)}\,d\mu(x)\,d\mu(y).
$$
We refer the reader to \cite{GK} and \cite{BS} for basic facts about Schatten--von Neumann classes.

\

{\bf Besov classes.} We consider here spaces ${\frak B}_p^s$ of distributions on $\R_+=(0,\be)$. 
The space ${\frak B}_p^s$ can be identified with the restrictions of the Fourier transforms of 
functions in the Besov classes $B_p^s(\R)$ to $\R_+$.

Let $v$ be a $C^\be$ function on $\R$ such that
$$
\supp v=\left[\frac12,2\right]\quad\mbox{and}\quad\sum_{j=-\be}^\be v\left(\frac{x}{2^j}\right)=1,\quad x>0.
$$
Put 
\bay
\label{vj}
v_j(x)=v\left(\frac{x}{2^j}\right),\quad x>0.
\ey

For $0<p\le\be$ and $s\in\R$, we define the space ${\frak B}_p^s$ as the space of distributions $\f$ on $\R_+$
such that 
$$
\|\f\|_{{\frak B}_p^s}\df\left(\sum_{j=-\be}^\be\left(2^{js}\|\F(v_j\f)\|_{L^p}\right)^p\right)^{1/p}<\be,\quad p<\be,
$$
and
$$
\|\f\|_{{\frak B}_\be^s}\df\sup_{-\be<j<\be}2^{js}\|\F(v_j\f)\|_{L^\be}<\be.
$$
Here $\F$ denotes the Fourier transform.

If $p\ge1$, ${\frak B}_p^s$ is a Banach space with norm $\|\cdot\|_{{\frak B}_p^s}$. If $p<1$, ${\frak B}_p^s$
is a quasi-Banach space with quasinorm $\|\cdot\|_{{\frak B}_p^s}$, i.e., $\|\f_1+\f_2\|_{{\frak B}_p^s}
\le\const\big(\|\f_1\|_{{\frak B}_p^s}+\|\f_2\|_{{\frak B}_p^s}\big)$.

We also define the space ${\frak b}_\be^s$, $s\in\R$, as the closed subspace of ${\frak B}_\be^s$, which consists 
of distributions $\f\in{\frak B}_\be^s$ such that
$$
\lim_{|j|\to\be}2^{js}\|\F(v_j\f)\|_{L^\be}=0.
$$

If $\f$ is a distribution on $\R_+$ and $\psi(x)=x^\s\f(x)$ 
(this equality has to be understood in the distributional sense), then it is easy to see that
$\f\in{\frak B}_\be^s$ if and only if $\psi\in{\frak B}_\be^{s-\s}$.

If $1\le p<\be$ and $s_1,\,s_2\in\R$, 
one can identify the dual space $\left({\frak B}_p^{s_1}\right)^*$ with the space ${\frak B}_{p'}^{s_2}$
with respect to the pairing
\bay
\label{bes}
\langle\f,\psi\rangle=
\int_0^\be t^{s_1+s_1}\f(t)\ov{\psi(t)}dt,\quad \f\in{\frak B}_p^{s_1},\quad\psi\in{\frak B}_{p'}^{s_2}.
\ey
Note that the integral on the right-hand side makes sense for compactly supported $C^\be$ functions $\f$ can 
be understood as the value of the distribution $\psi$ at the function $t\mapsto t^{s_1+s_1}\f(t)$. The linear
functional $\f\mapsto\langle\f,\psi\rangle$ extends by continuity to the whole space ${\frak B}_p^{s_1}$.
As usual, $\frac1{p}+\frac1{p'}=1$.

The dual space ${\frak b}_\be^{s_1}$ can be identified with the space ${\frak B}_1^{s_2}$ 
with respect to the same pairing \rf{bes}.

\medskip

{\bf Remark.} The Besov spaces $B_p^s(\R)$ of functions on $\R$ can be defined as the space of 
tempered distributions $f$ on $\R$ such that  
$$
\sum_{j=-\be}^\be \big(2^{sj}\|f*\xi_j\|_{L^p}\big)^p+\sum_{j=-\be}^\be \big(2^{sj}\|f*\eta_j\|_{L^p}\big)^p<\be,
$$
where $\xi_j$ and $\eta_j$ are functions in on $\R$ such that 
$$
\F\xi_j=v_j\quad\mbox{and}\quad(\F\eta_j)(x)=v_j(-x),\ \ x\in\R.
$$
and the $v_j$ are defined by \rf{vj}. This this space contains all polynomials.
It is possible to define the space $B_p^s(\R)$ modulo the polynomials of degree at most $s-1/p$, in which case 
only polynomials of degree at most $s-1/p$ can belong to this space. In both cases
$$
{\frak B}_p^s=\big\{\F f\big|(0,\be):~f\in B_p^s(\R)\big\}.
$$
The subspace $\big(B_p^s(\R)\big)_+$ of $B_p^s(\R)$ is defined by
$$
\big(B_p^s(\R)\big)_+=\big\{f\in B_p^s(\R):\supp\F f\subset[0,+\infty)\big\}.
$$
Functions in $\big(B_p^s(\R)\big)_+$ can be extended analytically in a natural way to the upper half-plane.
Clearly,
$$
\big\{(\F f)\big|(0,\be):~f\in\big(B_p^s(\R)\big)_+\big\}={\frak B}_p^s.
$$

\medskip

We refer the reader to \cite{Pe} for more information about Besov classes.

\

{\bf Weighted integral Hankel operators.} For a locally integrable 
function $\f$ on $\R_+$ the weighted integral Hankel operator
$\GG^{\a,\b}_\f$ is defined by 
$$
\big(\GG^{\a,\b}_\f f\big)(x)=\int_0^\be x^\a y^\b\f(x+y)f(y)dy
$$
for smooth functions $f$ with compact support in $\R_+$. Again, the definition makes sense
for distributions $\f$ on $\R_+$. The operators $\GG^{\a,\b}_\f$ are analogs of weighted
Hankel matrices $\G^{\a,\b}_\psi=\big\{(1+j)^\a(1+k)^\b\hat\psi(j+k)\big\}_{j,k\ge0}$, where $\psi$ is a function analytic in
the unit disk. 

For $\a=\b=0$, the operator $\GG^{0,0}_\f=\GG_\f$ is the integral Hankel operator defined in the introduction.

We need the following results.

\medskip

{\bf Theorem A.} {\it Suppose that $\a>0$ and $\b>0$. Then $\GG^{\a,\b}_\f$ is bounded on $L^2(\R_+)$
if and only if $\f\in{\frak B}_\be^{\a+\b}$ and $\GG^{\a,\b}_\f$ is compact if and only if $\f\in{\frak b}_\be^{\a+\b}$.}

\medskip

We refer the reader to \cite{P2} for the corresponding result for weighted Hankel matrices and to \cite{JP} for
more general results than Theorem A.

\medskip

{\bf Theorem B.} {\it Let $0<p<\be$. Suppose that $\min\{\a,\b\}>\max\{-\frac{1}{2},-\frac{1}{p}\}$. Then
$\GG^{\a,\b}_\f\in\bS_p$ if and only if $\f\in{\frak B}_p^{1/p+\a+\b}$.}

\medskip

For $p\ge1$ the description of the weighted Hankel matrices $\G^{\a,\b}_\psi$ of class \lb $\bS_p$ was
obtained in \cite{P1} for \mbox{$\a=\b=0$} and in \cite{P2} in the case 
\lb$\min\{\a,\b\}>\max\{-\frac{1}{2},-\frac{1}{p}\}$. For integral Hankel operators and $p\ge1$ see
\cite{CR}, \cite{R}, and \cite{JP}. In the case $p<1$ 
we refer the reader to \cite{P1} for weighted Hankel matrices and to
\cite{S} for the weighted integral Hankel operators. See also \cite{P}, \lb Ch. 6.

\

{\bf Averaging projection onto the integral Hankel operators.} Consider the orthogonal projection $\cp$ on 
the Hilbert--Schmidt class $\bS_2$ onto the subspace of Hankel integral operators. Clearly, if $k\in L^2(\R_+\times\R_+)$,
then 
$$
\cp{\frak I}_k=\GG_\f,
$$ 
where
$$
\f(x)=\frac1{x}\int_0^x k(t,x-t)dt.
$$

For $2<p<\be$, the averaging projection $\cp$ can be defined on the dense subset $\bS_p\cap\bS_2$ of of $\bS_p$.
As usual, we say that $\cp$ is bounded on $\bS_p$ if it extends by continuity to a bounded operator
on $\bS_p$.

\medskip

{\bf Theorem C.} {\it Let $1<p<\be$. Then the averaging projection $\cp$ is bounded on $\bS_p$}.

\medskip

The same result for the averaging projection onto the space of Hankel matrices was obtained in
\cite{P1}, see also \cite{P}. The same proof works for the averaging projection onto the set of integral Hankel operators.

Note that $\cp$ is unbounded on $\bS_1$ and on the spaces of bounded and compact operators (see \cite{P1}, \cite{P}).
However, there are bounded projections on $\bS_1$ onto the subspace of integral Hankel operators.
Indeed, for $\a,\,\b>0$ we define the weighted averaging projection $\cp_{\a,\b}$ by
$$
\cp_{\a,\b}{\frak I}_k=\GG_\f,
$$ 
where
$$
\f(x)=\frac{\int\limits_0^xt^\a (x-t)^\b k(t,x-t)dt}{\int\limits_0^xt^\a (x-t)^\b dt}.
$$

\medskip

{\bf Theorem D.} {\it Let $\a,\,\b>0$. Then $\cp_{\a,\b}$ is a bounded operator on $\bS_1$.}

\medskip

The corresponding result for weighted averaging projection on the subspace of Hankel matrices
was found by the first author, see
(\cite{P2} and \cite{P}). The same proof works in the case of integral operators.

\

\section{\bf Theorem B is sharp}
\setcounter{equation}{0}

\

When we study properties of the averaging projection $\Q_{\a,\b}$ onto the space of operators 
of the form ${\frak G}^{\a,\b}_\f$, we will need the fact that Theorem B cannot 
be extended to other values of $\a$ and $\b$.

\begin{thm}
\label{sharp}
Let $\a,\,\b\in\R$ and $0<p<\be$. Suppose that
$\min\{\a,\b\}\le\max\{-\frac{1}{2},-\frac{1}{p}\}$.
Then there are functions $\psi\in{\frak B}_p^{1/p+\a+\b}$ such that $\GG^{\a,\b}_\psi\not\in \bS_p$.
\end{thm}

\Pf Consider first the case $p\ge2$. For a positive integer $n$ we define the function $\f_n$ on $\R_+$
by 
$$
\f_n(x)=\left\{\begin{array}{ll}1,&x\in\left(1,1+\frac2{n}\right),\\0,&\mbox{otherwise}.\end{array}\right.
$$
It follows easily from the definition of ${\frak B}_p^{1/p+\a+\b}$ that
$$
\|\f_n\|_{{\frak B}_p^{1/p+\a+\b}}\le\const\|\F\f_n\|_{L^p}\le\const\cdot n^{-1/p'}.
$$

Let us now estimate from below $\big\|\GG^{\a,\b}_{\f_n}\big\|_{\bS_p}$. Define the function $k$ on
$\R_+^2$ by
$$
k(x,y)\df\left\{\begin{array}{ll}x^\a y^\b,&(x.y)\in\bigcup\limits_{k=0}^n\D_j,\\[.4cm]
0,&(x,y)\notin\bigcup\limits_{k=0}^n\D_j,
\end{array}\right.
$$
where
$$
\D_j=\left(\frac{j}{n},\frac{j+1}{n}\right)\times\left(\frac{n-j}{n},\frac{n-j+1}{n}\right),\quad 
0\le j\le n.
$$
Since the squares $\D_j$ have disjoint projections onto the coordinate axes, it is easy to see that
$$
\big\|\GG^{\a,\b}_{\f_n}\big\|_{\bS_p}
\ge\|{\frak I}_k\|_{\bS_p}=\left(\sum_{j=0}^n\|{\frak I}_{k_j}\|_{\bS_p}^p\right)^{1/p},
$$
where $k_j=k\chi_{\D_j}$ and $\chi_{\D_j}$ is the characteristic function of $\D_j$. 
Clearly, ${\frak I}_{k_j}$ is a rank one operator and
$$
\|{\frak I}_{k_j}\|_{\bS_p}\ge\frac{j^\a(n-j)^\b}{n^{\a+\b+1}}.
$$
Without loss of generality, we may assume $\a\le\b$. It is easy to verify that
\bey
\big\|\GG^{\a,\b}_{\f_n}\big\|_{\bS_p}
&\ge&\left(\sum_{j=0}^n\frac{j^{p\a}(n-j)^{p\b}}{n^{p(\a+\b+1)}}\right)^{1/p}\\[.2cm]
&\ge&\const\left\{\begin{array}{ll}n^{-\a-1},&\a<-1/p,\\[.1cm]n^{-\a-1}(\log(1+n))^{1/p},&\a=-1/p.
\end{array}\right.
\eey
Clearly,
$$
\frac{\big\|\GG^{\a,\b}_{\f_n}\big\|_{\bS_p}}{\|\f_n\|_{{\frak B}_p^{1/p+\a+\b}}}\to\be\quad\mbox{as}\quad n\to\be,
$$
which completes the proof in the case $p\ge2$.

Suppose now that $p<2$. Again, we assume that $\a\le\b$. We prove that if $\a\le-\frac12$, the condition
$\psi\in{\frak B}_p^{1/p+\a+\b}$ does not even imply that $\GG^{\a,\b}_\psi\in \bS_2$.
Let $\psi$ be a nonzero smooth function with support in $[1,2]$. Clearly, $\F\psi\in L^p$, and so
$\psi\in{\frak B}_p^{1/p+\a+\b}$. On the other hand,
$$
\big\|\GG^{\a,\b}_\psi\big\|^2_{\bS_2}=\int_1^2|\f(t)|^2\int_0^tx^{2\a}(t-x)^{2\b}dx\,dt=\be,
$$
since, clearly,
$$
\int_0^tx^{2\a}(t-x)^{2\b}dx=\be.\quad\bl
$$

\section{\bf Boundedness, Compactness, and Schatten--von Neumann Properties}
\setcounter{equation}{0}

\

Let $k$ be a function on $\R_+^2$ such that the integral operator on $L^2(\R_+)$ with kernel function
$k$ is bounded on $L^2(\R_+)$. As in \S 2, denote this integral operator by ${\frak I}_k$:
$$
({\frak I}_kf)(x)=\int_{\R_+}k(x,y)f(y)dy.
$$
We say that $k\in\bS_p(\R_+^2)$ if the operator ${\frak I}_k$ belongs to the Schatten--von Neumann class 
$\bS_p$ and we write $\|k\|_{\bS_p}\df\|{\frak I}_k\|_{\bS_p}$.
Let now $\a$ and $\b$ be nonzero real numbers. We put
\bay
\label{kab}
k_{\a,\b}(x,y)\df x^{\frac1{2\a}-\frac12}y^{\frac1{2\b}-\frac12}k\left(x^{\frac1{\a}},y^{\frac1{\b}}\right).
\ey
We introduce the unitary operator $U_\a$ on $L^2(\R_+)$ defined by
$$
(U_\a f)(x)=\frac1{\sqrt{|\a|}}x^{\frac1{2\a}-\frac12}f\left(x^{\frac1{\a}}\right),\quad f\in L^2(\R_+).
$$
It is easy to see that 
$$
U_\a {\frak I}_k=\frac1{\sqrt{|\a\b|}}{\frak I}_{k_{\a,\b}}U_\b,\quad\a,\,\b\in\R\setminus\{0\},
$$
and so
$$
\|k_{\a,\b}\|_{\bS_p}=\sqrt{|\a\b|}\cdot\|k\|_{\bS_p}.
$$

\begin{thm}
\label{bc}
Suppose that $\a,\,\b\in(0,1)$.
Then ${\frak G}^{\a,\b}_\f$ is bounded if and only if
$\f\in{\frak B}_\be^{\frac1{2\a}+\frac1{2\b}-1}$. The operator 
${\frak G}^{\a,\b}_\f$ is compact if and only if
$\f\in{\frak b}_\be^{\frac1{2\a}+\frac1{2\b}-1}$.
\end{thm}

\Pf  Suppose first that $\f$ is a locally integrable functin on $\R_+$.
Consider the kernel function $\vk$ of ${\frak G}^{\a,\b}$:
$$
\vk(x,y)=\f(x^\a+y^\b),\quad x,~y>0.
$$
Obviously, by \rf{kab}
\bay
\label{kappa}
\vk_{\a,\b}(x,y)=x^{\frac1{2\a}-\frac12}y^{\frac1{2\b}-\frac12}\f(x+y),
\ey
i.e., the integral operator operator with kernel function $\vk_{\a,\b}$ is a weighted integral
Hankel operator. Hence, ${\frak G}^{\a,\b}_\f$ is bounded (or compact) if and only if 
$\GG_\f^{\frac1{2\a}-\frac12,\frac1{2\b}-\frac12}$ 
is bounded (or compact). By Theorem A (see \S 2), this is equivalent to the fact that
$\f\in{\frak B}_\be^{\frac1{2\a}+\frac1{2\b}-1}$ $\Big($or $\f\in{\frak b}_\be^{\frac1{2\a}+\frac1{2\b}-1}\,\Big)$.

If $\f$ is a distribution, it is easy to verify that the formula
\bay
\label{dist}
\GG_\f^{\frac1{2\a}-\frac12,\frac1{2\b}-\frac12}=\sqrt{\a\b}\,U_\a{\frak G}^{\a,\b}_\f U^*_\b
\ey
still holds, which implies the result. $\bl$

\begin{thm}
\label{Sp}
Let $p$, $\a$, and $\b$ be positive numbers 
such that 
$$
\max\{\a,\b\}(p-2)<p.
$$
Then ${\frak G}^{\a,\b}_\f\in\bS_p$ if and only if
$\f\in{\frak B}_p^{\frac1{2\a}+\frac1{2\b}+\frac1{p}-1}$.
\end{thm}

\Pf Again, apply formula \rf{dist}.
Obviously, $\a$ and $\b$
satisfy the hypotheses of Theorem B. The result follows from Theorem B. $\bl$

It turns out that in the case $1\le p<\be$
the necessity of the condition \lb$\f\in{\frak B}_p^{\frac1{2\a}+\frac1{2\b}+\frac1{p}-1}$
holds for any positive $\a$ and $\b$.

\begin{thm}
\label{nc}
Let $1\le p<\be$, and let $\a$ and $\b$ be positive numbers. Suppose that 
${\frak G}^{\a,\b}_\f\in\bS_p$. Then $\f\in{\frak B}_p^{\frac1{2\a}+\frac1{2\b}+\frac1{p}-1}$.
\end{thm}

\Pf As we have already observed, the integral operator ${\frak I}_{\vk_{\a,\b}}$ with kernel function
$\vk_{\a,\b}$ defined by \rf{kappa} must belong to $\bS_p$.

Suppose first that $p>1$.
We apply to ${\frak I}_{\vk_{\a,\b}}$ the averaging projection $\cp$ (see \S 2).
It is easy to verify that 
$$
\cp{\frak I}_{\vk_{\a,\b}}=\const\GG_\psi,
$$ 
where 
\bay
\label{psi}
\psi(x)=x^{\frac1{2\a}+\frac1{2\b}-1}\f(x).
\ey
By Theorem C, $\GG_\psi\in\bS_p$. Now 
by Theorem B, $\psi\in{\frak B}_p^{1/p}$, which is equivalent to the fact that
$\f\in{\frak B}_p^{\frac1{2\a}+\frac1{2\b}+\frac1{p}-1}$ (see Section 2). 

Let now $p=1$. We apply to the operator ${\frak I}_{\vk_{\a,\b}}$ the 
weighted projection $\cp_{1.1}$, which is bounded on $\bS_1$ (see Theorem D). It is easy to verify that 
$$
\cp_{1,1}{\frak I}_{\vk_{\a,\b}}=\const\GG_\psi,
$$ 
where $\psi$ is defined by \rf{psi}. By Theorem B, $\psi\in{\frak B}_1^1$, which is equivalent 
to the fact that $\f\in{\frak B}_1^{\frac1{2\a}+\frac1{2\b}}$ (see \S 2). $\bl$

\

\section{\bf The Averaging Projection onto the Operators $\bs{{\frak G}^{\a,\b}_\f}$}
\setcounter{equation}{0}

\

In this section we study metric properties of the averaging projection on the class of 
operators of the form ${\frak G}^{\a,\b}_\f$. Consider the orthogonal projection
$\Q_{\a,\b}$ from the Hilbert--Schmidt class $\bS_2$ onto the subspace of $\bS_2$
of operators of the form ${\frak G}^{\a,\b}_\f$.
Clearly, $\Q_{1,1}$ is just the averaging projection $\cp$ onto the Hankel integral operators.
If we identify the Hilbert--Schmidt operators on $L^2(\R_+)$ with the space $L^2(\R_+^2)$, 
we find that $\Q_{\a,\b}$ is the orthogonal projection onto the subspace of functions
that are constant on the sets $\big\{(x,y)\in\R_+^2:~x^\a+y^\b=c\big\}$, $c>0$.

We are going to characterize those $\a,\,\b$, and $p$, for which the projection
$\Q_{\a,\b}$ is bounded on $\bS_p$.
Clearly, $\Q_{\a,\b}T$ is well-defined for $T\in\bS_p$ if $p\le2$. If $p>2$, we say that $\Q_{\a,\b}$
is bounded on $\bS_p$ if it extends to a bounded operator from $\bS_p\cap\bS_2$.
Since $\Q_{\a,\b}$ is self-adjoint on $\bS_2$, it follows that $\Q_{\a,\b}$ is bounded on $\bS_p$ if and only if 
it is bounded on $\bS_{p'}$, $1<p<\be$, and $\Q_{\a,\b}$ is bounded on $\bS_1$ if and only if it is
bounded on $\B(L^2(\R_+))$.
 
\begin{lem}
\label{lap}
Let $k\in L^2(\R_+^2)$ and let $\a,\,\b>0$. Put
\bay
\label{sym}
\f(x)\df\frac{\int\limits_0^{\frac{\pi}2}k\big(x^{\frac1{\a}}\cos^{\frac2{\a}}t,x^{\frac1{\b}}\sin^{\frac2{\b}}t\big)
\cos^{\frac2{\a}-1}t\sin^{\frac2{\b}-1}t\,dt}
{A(\a,\b)},
\ey
where
$$
A(\a,\b)\df
\int\limits_0^{\frac{\pi}2}\cos^{\frac2{\a}-1}t\sin^{\frac2{\b}-1}t\,dt
$$
Then 
\bay
\label{dap}
\Q_{\a,\b}{\frak I}_k={\frak G}^{\a,\b}_\f.
\ey
\end{lem}

\Pf The result follows from the following easily verifiable formula:
\bay
\label{fxy}
\iint\limits_{\R_+^2}f(x,y)dx\,dy=
\frac2{\a\b}\int\limits_0^\be \!r^{\frac1{\a}+\frac1{\b}-1}\!\!\int\limits_0^{\frac{\pi}2}
\!\!f\big(r^{\frac1{\a}}\cos^{\frac2{\a}}t,r^{\frac1{\b}}\sin^{\frac2{\b}}t\big)
\cos^{\frac2{\a}-1}t\sin^{\frac2{\b}-1}t\,dt\,dr,\nonumber\\
\ey
which holds for any nonnegative measurable function $f$. $\bl$

\medskip

{\bf Remark.} Note that $A(\a,\b)=\frac12B(1/\alpha,1/\beta)$, where $B$ is the Euler Beta function. 

\medskip

\begin{thm}
\label{ap}
Let $1\le p<\be$, and let $\a$ and $\b$ be positive numbers. Suppose that
\bay
\label{ogr}
-p<\max\{\a,\b\}(p-2)<p.
\ey
Then $\Q_{\a,\b}$ is bounded on $\bS_p$.
\end{thm}

Note that if $\max\{\a,\b\}\le1$, the rightmost inequality in \rf{ogr} holds for any $p$. If 
$\max\{\a,\b\}>1$, then \rf{ogr} is equivalent to the inequalities 
$$
\frac{2\max\{\a,\b\}}{\max\{\a,\b\}+1}<p<\frac{2\max\{\a,\b\}}{\max\{\a,\b\}-1}.
$$
Clearly, $\frac{2\max\{\a,\b\}}{\max\{\a,\b\}+1}$ and $\frac{2\max\{\a,\b\}}{\max\{\a,\b\}-1}$ are 
dual exponents, i.e., the sum of their reciprocals is euqal to one.

\Pf Consider first the case $p>1$. 
Let ${\mathcal T}_{\a,\b}$ be the operator on $\bS_p$ defined by
$$
{\mathcal T}_{\a,\b}T=\f,\quad\mbox{where}\quad \Q_{\a,\b}T={\frak G}^{\a,\b}_\f,
$$
see Lemma \ref{lap}. By Theorem \ref{Sp}, we have to show that 
${\mathcal T}_{\a,\b}$ is a bounded operator from  $\bS_p$ to ${\frak B}_p^{\frac1{p}+\frac1{2\a}+\frac1{2\b}-1}$.

Consider the dual exponent $p'=p/(p-1)$. Define the operator 
$$
\cR:{\frak B}_{p'}^{\frac1{p'}+\frac1{2\a}+\frac1{2\b}-1}\to\bS_{p'}(\R_+^2)
$$
by
\bay
\label{R}
(\cR\f)(x,y)=\f(x^\a+y^\b),\quad x,~y>0.
\ey
By Theorem \ref{Sp}, $\cR$ is a bounded operator. %and
%$$
%\|\cR\f\|_{\bS_{p'}}\ge\const\|\f\|_{{\frak B}_{p'}^{\frac1{p'}+\frac1{2\a}+\frac1{2\b}-1}}.
%$$

Let $k\in\bS_1(\R_+^2)$ and let
$$
\Q_{\a,\b}{\frak I}_k={\frak G}^{\a,\b}_\psi
$$
(see \rf{sym} and \rf{fxy}).
We have
\bay
\label{dual}
\iint\limits_{\R_+^2} (\cR\f)(x,y)\ov{k(x,y)}dx\,dy&=&
\iint\limits_{\R_+^2}(\cR\f)(x,y)\ov{\psi(x^\a+y^\b)}dx\,dy\nonumber\\[.2cm]
&=&\frac{2A(\a,\b)}{\a\b}\int\limits_0^\be r^{\frac1{\a}+\frac1{\b}-1}\f(r)\ov{\psi(r)}dr.
\ey

It follows that ${\mathcal T}_{\a,\b}=\cR^*$ if we identify 
$\left({\frak B}_{p'}^{\frac1{p'}+\frac1{2\a}+\frac1{2\b}-1}\right)^*$ with
${\frak B}_{p}^{\frac1{p}+\frac1{2\a}+\frac1{2\b}-1}$ with respect to the pairing
\bay
\label{pair}
\langle f,g\rangle=\int_{\R_+}t^{\frac1{\a}+\frac1{\b}-1}f(r)\ov{g(r)}dr
\ey
(see \S 2).
Hence, ${\mathcal T}_{\a,\b}$ is a bounded.

Suppose now that $p=1$. Clearly, \rf{ogr} means that $0<\a,\,\b<1$. Consider the operator 
$$
\cR:{\frak b}_{\be}^{\frac1{2\a}+\frac1{2\b}-1}\to\bs{\cc}(L^2(\R_+))
$$
defined by \rf{R}. It is bounded by Theorem A. 
Again, it is easy to see from \rf{dual} that ${\mathcal T}_{\a,\b}^*=\cR$ if we identify 
$\left({\frak b}_{\be}^{\frac1{2\a}+\frac1{2\b}-1}\right)^*$ with
${\frak B}_{1}^{\frac1{2\a}+\frac1{2\b}}$ with respect to the pairing \rf{pair}.
$\bl$

Now suppose that $\a,\,\b\in(0,1)$. By Theorem \ref{ap}, the averaging projection $\Q_{\a,\b}$ is bounded
on $\bS_1$. We can consider the conjugate operator $\Q^*_{\a,\b}$ on the space $\B(L^2(\R_+))$
with respect to the standard pairing between $\bS_1$ and $\B(L^2(\R_+))$ (see \S 2). Since $\Q_{\a,\b}$ is a self-adjoint operator
on $\bS_2$, it follows that $\Q^*_{\a,\b}T=\Q_{\a,\b}T$ for $T\in\bS_2$. We can now extend the averaging projection
$\Q_{\a,\b}$ to the space $\B(L^2(\R_+))$ by 
$$
\Q_{\a,\b}T=\Q^*_{\a,\b}T,\quad T\in \B(L^2(\R_+)).
$$
It is easy to show that if $T$ is a bounded integral operator on $L^2(\R_+)$ with kernel function $k$, then
$\Q_{\a,\b}T$ can be defined as in \rf{sym} and \rf{dap}.

\begin{thm}
\label{be}
Let $\a,\,\b\in(0,1)$. Then $\Q_{\a,\b}$ is a bounded operator on the space $\B(L^2(\R_+))$ 
of bounded operators and on the space $\bs{\cc}(L^2(\R_+))$ of compact operators.
\end{thm}

\Pf We have already explained the fact that $\Q_{\a,\b}$ is bounded on $\B(L^2(\R_+))$.
The boundedness of $\Q_{\a,\b}$ on $\bs{\cc}(L^2(\R_+))$  follows immediately from the fact
that $\Q_{\a,\b}\bS_2\subset\bS_2$ and the fact that $\bS_2$ is dense in $\bs{\cc}(L^2(\R_+))$. $\bl$

Let us prove now that Theorems \ref{ap} and \ref{be} are sharp.

\begin{thm}
\label{0pbe}
Suppose that $\a$ and $\b$ are positive numbers such that the averaging projection 
$\Q_{\a,\b}$ is bounded on $\bS_p$, $0<p<\be$. Then $p\ge1$ and \rf{ogr} holds.
\end{thm}

\begin{thm}
\label{pbe}
Suppose that $\a$ and $\b$ are positive numbers such that the averaging projection 
$\Q_{\a,\b}$ is bounded on the space $\B(L^2(\R_+))$ or on the space $\bs{\cc}(L^2(\R_+))$.
then $\a,\,\b<1$.
\end{thm}

In fact for $p<1$ the following much stronger retsult holds.

\begin{thm}
\label{p<1}
Let $0<p<1$ and $\a,\,\b>0$. Then there is no bounded projection 
from $\bS_p$ onto the subspace of operators of the form
${\frak G}^{\a,\b}_\f$.
\end{thm}

Let us first prove Theorem \ref{p<1}.

\medskip

{\bf Proof of Theorem \ref{p<1}.} The result follows from Theorem \ref{Sp} and the Kalton
theorem \cite{K}, which says, in particular, that if $X$ is a coplemented subspace of $\bS_p$, $0<p<1$, such that
$X$ can be imbedded isomorphically to an $L^p$ space, then the $\bS_p$ quasinorm and the $\bS_1$ norm 
on $X$ are equivalent. Indeed, let $X$ be the subspace of $\bS_p$ of operators of the form ${\frak G}^{\a,\b}_\f$.
By Theorem \ref{Sp}, the $\bS_1$ norm and the $\bS_p$ quasinorm on $X$ are not equivalent.
It follows easily from the definition of the spaces ${\frak B}_p^s$ given
in \S 2 that $X$ can be imbedded isometrically in an $L^p$ space. Thus by the Kalton theorem
$X$ is not a complemented subspace of $\bS_p$. $\bl$

\medskip

{\bf Proof of Theorem \ref{0pbe}.} Suppose that $1\le p<\be$.
The reasoning given in the proof of Theorem \ref{ap} shows that if the averaging projection 
$\Q_{\a,\b}$ is bounded on $\bS_p$, then the condition $\f\in{\frak B}_p^{\frac1{2\a}+\frac1{2\b}+\frac1{p}-1}$
implies ${\frak G}^{\a,\b}_\f\in\bS_p$. The result now follows from \rf{kappa} and from Theorem
\ref{sharp}. $\bl$

\medskip

{\bf Proof of Theorem \ref{pbe}.} If $\Q_{\a,\b}$  is bounded on 
$\B(L^2(\R_+))$ or on $\bs{\cc}(L^2(\R_+))$, then, by duality, $\Q_{\a,\b}$
is bounded on $\bS_1$. The result follows now from Theorem \ref{0pbe}. $\bl$

\

\

\noindent
\begin{tabular}{p{8cm}p{14cm}}
A.B. Aleksandrov & V.V. Peller \\
St-Petersburg Branch & Department of Mathematics \\
Steklov Institute of Mathematics  & Michigan State University \\
Fontanka 27, 191011 St-Petersburg & East Lansing, Michigan 66506\\
Russia&USA%\\
%&and\\
%&Equipe d'Analyse\\
%&Universit\'{e} Paris VI, Boite 186\\
%&75252 Paris Cedex 05\\
%&France
\end{tabular}


\begin{thebibliography}{99}

\bibitem[AJPR]{AJPR} {\sc A.B. Aleksandrov, S. Janson, V.V. Peller,} and {\sc R. Rochberg},
{\em An interesting class of operators with unusual Schatten--von Neumann behavior}. In: 
Function Spaces, Interpolation Theory and Related Topics, Proceedings of the International 
Conference in Honour of Jaak Peetre on his 65th Birthday, 61--149, Walter de Gruyter,
Berlin, 2002.

\bibitem[BS]{BS} {\sc M.S. Birman} and {\sc M.Z. Solomyak},
{\em Spectral theory
of self-adjoint operators in Hilbert space}, Reidel Publishing Company,
Dordrecht, 1986.

\bibitem[CR]{CR}
{\sc R.R. Coifman} and {\sc R. Rochberg}, 
{\it Representation theorems 
for holomorphic and harmonic functions in $L^p$, 
representation theorems for Hardy spaces,} Ast\'{e}risque {\bf 77} (1980), 11--66.

\bibitem[GK]{GK} {\sc I.C. Gohberg and M.G. Krein}, {\it Introduction to the theory
of linear nonselfadjoint operators in Hilbert space,} Nauka, Moscow, 1965; 
English transl.: Amer. Math. Soc., Providence, RI, 1969.

\bibitem[JP]{JP} {\sc S. Janson} and {\sc J. Peetre},
{\it Paracommutators -- boundedness
and Schatten--von Neumann classes,} Trans. Amer. Math. Soc. {\bf 305}
(1988), 467--504.

\bibitem[K]{K} {\sc N.J. Kalton}, {\em Plurisubharmonic functions on
quasi-Banach spaces}, Studia Math. {\bf84} (1986), no. 3, 297--324.

\bibitem[Pee]{Pe} {\sc J. Peetre}, 
{\em New thoughts on Besov spaces}, Duke Univ. Press., Durham, NC, 1976.

\bibitem[Pel1]{P1} {\sc V.V. Peller}, {\it Hankel operators of class ${\frak S}_p$
and applications (rational approximation, Gaussian processes, the majorization
problem for operators),} Mat Sb. {\bf 113} (1980), 538-581.
English transl.: Math. USSR-Sb. {\bf 41} (1982), 443-479.

\bibitem[Pel2]{P2} {\sc V.V. Peller},
{\em Vectorial Hankel operators and related operators
of the Schatten--von Neumann class ${\bf S}_{p}$}, Int. Equat. Oper. 
Theory {\bf 5} (1982), 244--272.
 
\bibitem[Pel3]{P3} {\sc V.V. Peller}, {\it A description of Hankel operators of class
${\frak S}_p$ for $p>0$, investigation of the rate of rational approximation
and other applications,} Mat. Sb. {\bf 122} (1983), 481-510.
English transl.: Math. USSR-Sb. {\bf 50} (1985), 465-494.

\bibitem[Pel4]{P} {\sc V.V. Peller}, {\em Hankel operators and their applications},
to appear in Springer-Verlag.

\bibitem[R]{R} {\sc R. Rochberg},
{\it Trace ideal criteria for Hankel operators and
commutators,} Indiana Univ. Math. J. {\bf31} (1982), 913--925.

\bibitem[S]{S}{\sc S. Semmes},
{\it Trace ideal criteria for Hankel operators and applications to Besov 
spaces,}  Integral Equations and Operator Theory {\bf 7} (1984), 241--281.

\end{thebibliography}
\end{document}